\documentclass[10pt]{article}

\usepackage{amssymb,amsfonts}

\newcommand{\CC}{\mathbb{C}}

\newcommand{\RR}{\mathbb{R}}
\newcommand{\Sn}{\mathcal{S}_{n}}
\newcommand{\SE}{\mathbb{S}_{n}}
\newcommand{\qed}{\hbox{}\nobreak\hfill\vrule width 1.4mm height 1.4mm depth 0mm
    \par \goodbreak \smallskip}

\newtheorem{theorem}{Theorem}
\newtheorem{proposition}[theorem]{Proposition}

\newtheorem{definition}{Definition}
\newenvironment{proof}{\noindent{\it Proof.} \rm }{\hfill\qed}

\begin{document}

\title{On Reidys and Stadler's metrics for RNA secondary 
structures\thanks{This work has been partially supported by the 
Spanish DGES, grant BFM2000-1113-C02-01.}}

\author{F. Rossell\'o\\[1ex]
{\small Departament de Matem\`atiques i Inform\`atica,}\\
{\small Institut Universitari d'Investigaci\'o en Ci\`encies de la 
Salut (IUNICS),}\\
{\small Universitat de les Illes Balears,}\\
{\small 07122 Palma de Mallorca (Spain)}\\
{\small \emph{E-mail:} \texttt{cesc.rossello@uib.es}}}

\date{}

\maketitle

\begin{abstract}
We compute explicitly several abstract metrics  for RNA secondary 
structures defined by Reidys and Stadler.
\smallskip

\noindent\textbf{Keywords:} RNA secondary structure, metric, 
symmetric group.
\end{abstract}

\section{Introduction}

As it is well known, an RNA molecule can be viewed as a chain of 
(ribo)nucleotides with a definite orientation.  Each of these 
nucleotides is characterized by (and in practice identified with) the 
base attached to it, which can be adenine (A), cytosine (C), guanine 
(G), or uracil (U).  Thus, an RNA molecule with $N$ nucleotides can be 
mathematically described as a word of length $N$ over the alphabet 
$\{A,C,G,U\}$, called the \emph{primary structure} of the molecule.

In the cell and \textsl{in vitro} each RNA molecule folds into a 
three-dimensional structure, which determines its biochemical 
function.  This  structure is held together by weak 
interactions called \emph{hydrogen bonds} between pairs of 
non-consecutive bases: actually, a hydrogen bond can only form between 
bases that are several positions apart in the chain, but we shall 
not take this restriction into account here.  Most of these bonds form 
between \emph{Watson-Crick complementary bases}, i.e., between $A$ and 
$U$ and between $C$ and $G$, but a significant amount of bonds also 
form between other pairs of bases \cite{beyondWC}.  The 
\emph{secondary structure} of an RNA molecule is a simplified model of 
this three-dimensional structure, consisting of an undirected graph 
with nodes its bases and arcs its base pairs or \emph{contacts}; the 
\emph{length} of a secondary structure is the number of its 
nodes.  A restriction is added to the definition of secondary 
structure: a base can only pair with at most one base.  This restriction is 
called the \emph{unique bonds condition}.

An important problem in molecular biology is the comparison of these 
RNA secondary structures, because it is assumed that a preserved 
three-dimensional structure corresponds to a preserved function.  
Moreover, the comparison of RNA secondary structures of a fixed length is 
used in the prediction of RNA secondary structures to reduce the 
output of alternate structures when suboptimal solutions, and not 
only optimal, are considered \cite[\S IX]{Zuk}.  In a seminal paper on 
the algebraic representation of biomolecular structures \cite{RS96}, 
C. Reidys and P. F. Stadler introduced three abstract metrics on the 
set of RNA secondary structures of a fixed length based on their 
algebraic models and independent of any notion of graph edition, and 
they discussed their biophysical relevance.  They ended that paper by 
asking, among other questions, whether there exists any relation 
between the metrics for RNA secondary structures they had defined.  In 
this paper we answer this question by explicitly computing these 
metrics.  In a subsequent paper \cite{LR} we plan to generalize these 
metrics to contact structures without unique bonds, as for instance 
protein structures.

\section{Main results}

From now on, let $[n]$ denote the set $\{1,\ldots,n\}$, for every 
positive integer $n$.

\begin{definition}
An \emph{RNA secondary structure} of length $n$ is an undirected graph 
without multiple edges or self-loops $\Gamma=([n],Q)$, for some $n\geq 
1$, whose arcs $\{j,k\}\in Q$, called \emph{contacts}, satisfy the 
following two conditions:

i)  For every $j\in [n]$, $\{j,j+1\}\notin Q$.

ii) For every $j\in [n]$, if $\{j,k\},\{j,l\}\in Q$, then $k=l$.
\end{definition}

Condition (i) translates the impossibility of a contact between two 
consecutive bases, while condition (ii) translates the {unique bonds 
condition}.  We should point out that this definition of RNA secondary 
structure is not the usual one, as the latter forbids the existence of 
(\emph{pseudo})\emph{knots}: pairs of contacts $\{i,j\}$ and $\{k,l\}$ 
such that $i<k<j<l$.  This rather unnatural condition is usually 
required in order to enable the use of dynamic programming methods to 
predict RNA secondary structures \cite{Zuk}, but real secondary 
structures can contain knots and thus we shall not impose this 
restriction here.  Therefore, our RNA secondary structures correspond 
to what in the literature on secondary structure modelling has been 
called \emph{contact structures with unique bonds} \cite{RS96,SS99} or 
\emph{1-diagrams} \cite{HS99}.

We shall denote from now on a contact $\{j,k\}$ by $j\!\cdot\!  k$ or 
$k\!\cdot\!  j$, without distinction.  A node is said to be 
\emph{isolated} in an RNA secondary structure when it is not involved 
in any contact.

Let $\SE$ stand for the set of all RNA secondary structures of length 
$n$ and let $\Sn$ be the symmetric group of permutations of $[n]$.

\begin{definition}
For every $\Gamma=([n],Q)\in \SE$,  say with 
$Q=\{i_{1}\!\cdot\!j_{1},\ldots,i_{k}\!\cdot \!j_{k}\}$, let
$$
\pi(\Gamma)=\prod_{t=1}^k (i_{t},j_{t})\in \Sn,
$$
where $(i,j)$ denotes the transposition in $\Sn$ defined by $i\leftrightarrow j$.
\end{definition}

Reidys and Stadler proved in \cite{RS96} that the mapping $\pi:\SE\to 
\Sn$ is injective and that $\pi(\Gamma)$ is an involution for every 
$\Gamma\in \SE$.  This representation of RNA secondary structures as 
involutions is then used by these authors to define the following 
metric, called the \emph{involution metric}.

\begin{proposition}\label{inv-dist}
The mapping $d_{inv}:\SE\times \SE\to \RR$ sending every 
$(\Gamma_{1},\Gamma_{2})\in \SE^2$ to the least number 
$d_{inv}(\Gamma_{1},\Gamma_{2}) $ of transpositions necessary to represent 
the permutation $\pi(\Gamma_{1})\pi(\Gamma_{2})$, is a metric.
\end{proposition}

The following proposition computes explicitly 
this metric.  In it, and henceforth, $A\Delta B$ denotes the 
symmetric difference $(A\cup B)-(A\cap B)$ of the sets $A$ and $B$, and 
$|A|$ stands for the cardinal of the finite set $A$.

\begin{proposition}
For every $\Gamma_{1}=([n],Q_{1}),\Gamma_{2}=([n],Q_{2})\in \SE$,
$$
d_{inv}(\Gamma_{1},\Gamma_{2})=|Q_{1}\Delta Q_{2}|-2\Omega,
$$
where $\Omega$ is the number of cyclic orbits of length greater than 2 
induced by the action on $[n]$ of the subgroup $\langle 
\pi(\Gamma_{1}),\pi(\Gamma_{2})\rangle$ of $\Sn$.
\end{proposition}

\begin{proof}
Let $\Gamma_{1}=([n],Q_{1})$ and $\Gamma_{2}=([n],Q_{2})$ be two RNA 
secondary structures of length $n$.  To simplify the language, we 
shall refer to the orbits induced by the action of $\langle 
\pi(\Gamma_{1}),\pi(\Gamma_{2})\rangle$ on $[n]$ simply by 
\emph{orbits}. Notice that we can understand such an orbit as a subset 
$\{i_{1},i_{2},\ldots,i_{m}\}$ of $[n]$,
$m\geq 1$, such that
$$
i_{1}\!\cdot\!  i_{2},i_{2}\!\cdot\!  i_{3},\ldots,
i_{m-1}\!\cdot\!i_{m} \in Q_{1}\cup Q_{2}
$$
and maximal with this property, i.e., such that any other contact
in $Q_{1}\cup Q_{2}$ involving $i_{1}$ or $i_{m}$ can only be
$i_{1}\!\cdot \!  i_{m}$. The unique bonds condition (or, in 
group-theoretical terms, the fact that the transpositions defining 
each $\pi(\Gamma_{i})$ are pairwise disjoint) implies that 
if $\{i_{1},i_{2},\ldots,i_{m}\}$ is an orbit, then either
$$
i_{1}\!\cdot\!
i_{2},i_{3}\!\cdot\!
i_{4},\ldots,\in
Q_{1}\mbox{ and }
i_{2}\!\cdot\!i_{3},i_{4}\!\cdot\!i_{5},\ldots\in
Q_{2}
$$
or
$$
i_{1}\!\cdot\!
i_{2},i_{3}\!\cdot\!
i_{4},\ldots,\in
Q_{2}\mbox{ and }
i_{2}\!\cdot\!i_{3},i_{4}\!\cdot\!i_{5},\ldots\in
Q_{1}.
$$
Such an orbit is \emph{cyclic} if $m=2$ and $i_{1}\!\cdot\!  i_{2}\in 
Q_{1}\cap Q_{2}$, or $m\geq 3$ and $i_{1}\!\cdot\!  i_{m}\in Q_{1}\cup 
Q_{2}$, and an orbit is \emph{linear} in all other cases.  The fact 
that $\pi(\Gamma_{1}),\pi(\Gamma_{2})$ are both involutions implies 
that the cardinal of cyclic orbits is always even: roughly speaking, 
if $i_{1}\!\cdot\!  i_{2}\in Q_{1}$ in a cyclic orbit, then 
$i_{1}\!\cdot\!  i_{m}\in Q_{2}$ and hence $i_{m-1}\!\cdot\!  i_{m}\in 
Q_{1}$.

If two transpositions appearing in the product 
$\pi(\Gamma_{1})\pi(\Gamma_{2})$ are not disjoint, then the indexes 
involved in them belong to the same orbit.  Moreover, two 
disjoint transpositions always commute.  This allows us to reorganize 
the transpositions in the product $\pi(\Gamma_{1})\pi(\Gamma_{2})$, 
assembling them into subproducts corresponding to orbits.  More 
specifically, 
if for every orbit $O$ and for every $i=1,2$ we let
$$
\pi(O,\Gamma_{i})=\prod_{\scriptstyle k\cdot l\in Q_{i}\atop \scriptstyle k,l\in O} 
(k,l),\
$$
then
$$
\pi(\Gamma_{1})\pi(\Gamma_{2})=\prod_{O\in\{\mathrm{orbits}\}}
\pi(O,\Gamma_{1})\pi(O,\Gamma_{2}).
$$
Since the orbits are pairwise disjoint, this finally shows that the 
least number of transpositions which $\pi(\Gamma_{1})\pi(\Gamma_{2})$ 
decomposes into is equal to the sum of the least numbers of 
transpositions which $\pi(O,\Gamma_{1})\pi(O,\Gamma_{2})$ decompose 
into, for every orbit $O$.  It remains to compute this last number for 
each type of orbit $O$.

If $O$ is a linear orbit of length $m=1$, then 
$\pi(O,\Gamma_{1})\pi(O,\Gamma_{2})=\mathrm{Id}$, and it corresponds 
to a node that is isolated both in $\Gamma_{1}$ and in $\Gamma_{2}$.

Let now $O=\{i_{1},\ldots,i_{m}\}$ be a linear orbit of length $m\geq 
2$.  Consider first the case when $i_{1}\!\cdot\!  
i_{2},i_{3}\!\cdot\!  i_{4},\ldots,i_{m-1}\!\cdot\!i_{m}\in Q_{1}$ and 
$i_{2}\!\cdot\!i_{3},i_{4}\!\cdot\!i_{5},\ldots\in Q_{2}$; in 
particular, $m$ is even.  Then
$$
\begin{array}{rl}
\pi(O,\Gamma_{1})\pi(O,\Gamma_{2}) & = (i_{1},i_{2})(i_{3},i_{4})\cdots 
(i_{m-1},i_{m})(i_{2},i_{3})\cdots (i_{m-2},i_{m-1})\\ & =
(i_{2},i_{4},\ldots,i_{m},i_{m-1},i_{m-3},\ldots,i_{3},i_{1}),
\end{array}
$$
a cycle of length $m$ that decomposes into the product of $m-1$ 
transpositions (and it is the least number of transpositions required 
to represent it), which is exactly the number of contacts of 
$Q_{1}\cup Q_{2}$ involved in this orbit.

A similar argument shows that in all other cases for a linear orbit $O$, 
the permutation $\pi(O,\Gamma_{1})\pi(O,\Gamma_{2})$ is equal to a 
cycle of length the number of elements of the orbit, and thus the 
least number of transpositions this product decomposes into is equal to 
the number of contacts of $Q_{1}\cup Q_{2}$ involved in this orbit 
$O$, all of them belonging to $Q_{1}\Delta Q_{2}$.

If $O$ is a cyclic orbit of length $m=2$, say $O=\{i_{1},i_{2}\}$, 
then $\pi(O,\Gamma_{1})\pi(O,\Gamma_{2})=(i_{1},i_{2})(i_{1},i_{2})= 
\mathrm{Id}$.  Notice that cyclic orbits of length 2 correspond to 
contacts in $Q_{1}\cap Q_{2}$.

Finally, assume that $O$ is a cyclic orbit of length $m\geq 3$, say 
$O=\{i_{1},\ldots,i_{m}\}$ with $i_{1}\!\cdot\!  i_{2},i_{3}\!\cdot\!  
i_{4},\ldots,i_{m-1}\!\cdot\!i_{m}\in Q_{1}$ and 
$i_{2}\!\cdot\!i_{3},\ldots,i_{m-2}\!\cdot\!  i_{m-1},i_{m}\!\cdot\!  
i_{1}\in Q_{2}$; remember that $m$ is in this case even.
Then
$$
\begin{array}{rl}
\pi(O,\Gamma_{1})\pi(O,\Gamma_{2}) & = (i_{1},i_{2})(i_{3},i_{4})\cdots 
(i_{m-1},i_{m})(i_{2},i_{3})\cdots (i_{m-2},i_{m-1})(i_{m},i_{1})\\ & =
(i_{2},i_{4},\ldots,i_{m})(i_{m-1},i_{m-3},\ldots,i_{3},i_{1}),
\end{array}
$$
the product of two disjoint cycles of length $m/2$.  Since each cycle 
requires $m/2-1$ transpositions, the least number of transpositions 
the permutation $\pi(O,\Gamma_{1})\pi(O,\Gamma_{2})$ decomposes into 
is equal to $m-2$, the number of contacts of $Q_{1}\cup Q_{2}$ 
involved in this orbit $O$ (all of them belonging again to 
$Q_{1}\Delta Q_{2}$) minus 2.

To sum up, and if we call $\Omega$ the number of cyclic orbits of 
length greater than 2,
$$
\begin{array}{rl}
d_{inv}(\Gamma_{1},\Gamma_{2}) & =|\{\mbox{contacts involved in linear orbits}\}|\\
& \qquad +|\{\mbox{contacts involved in cyclic orbits 
of length greater than 2}\}| -2\Omega \\
& =|Q_{1}\Delta Q_{2}|-2\Omega,
\end{array}
$$
as we claimed.
\end{proof}

The number and structure of the orbits induced by the action of $\langle 
\pi(\Gamma_{1}),\pi(\Gamma_{2})\rangle$ on $[n]$ are related to the 
probability of transition from the neutral network of $\Gamma_{1}$ 
(the set of sequences that fold into it) to that of $\Gamma_{2}$: see 
\cite[\S 3]{RS96} and the references cited therein.

Let now $\mbox{Sub}(\Sn)$ be the set of subgroups of 
$\Sn$.

\begin{definition}
For every $\Gamma=([n],Q)\in \SE$, say with 
$Q=\{i_{1}\!\cdot\!j_{1},\ldots,i_{k}\!\cdot \!j_{k}\}$,
let
$$
T(\Gamma)=\{(i_{1},j_{1}),\ldots,(i_{k},j_{k})\}
$$
be the set of the transpositions corresponding to the contacts in $Q$ 
and let $G(\Gamma)=\langle T(\Gamma)\rangle$ be 
the subgroup of $\Sn$ generated by this set of transpositions.
\end{definition}

Reidys and Stadler also proved in \cite{RS96} that the mapping 
$G:\SE\to \mbox{Sub}(\Sn)$ is injective, and then they used this 
representation of RNA secondary structures as permutation subgroups  to define the following \emph{subgroup 
metric}.

\begin{proposition}
The mapping $d_{sgr}:  \SE\times \SE  \to  \RR$ defined by
$$
d_{sgr}(\Gamma_{1},\Gamma_{2})= \ln\left(\frac{\displaystyle |G(\Gamma_{1})\cdot 
G(\Gamma_{2})|}{\displaystyle |G(\Gamma_{1})\cap G(\Gamma_{2})|}\right)
$$
is a metric.
\end{proposition}

Next proposition shows that this metric simply measures, up to a 
constant factor, the cardinal of the symmetric difference of the sets 
of contacts.

\begin{proposition}
For every $\Gamma_{1}=([n],Q_{1}),\Gamma_{2}=([n],Q_{2})\in \SE$,
$$
d_{sgr}(\Gamma_{1},\Gamma_{2})= (\ln 2)|Q_{1}\Delta Q_{2}|.
$$
\end{proposition}

\begin{proof}
Since the transpositions generating a group $G(\Gamma)$, with 
$\Gamma\in \SE$, are pairwise disjoint, there is a bijection between 
$G(\Gamma)$ and the powerset $\mathcal{P}(T(\Gamma))$: each element of 
$G(\Gamma)$ is the product of a subset of $T(\Gamma)$ in a unique way.
Hence, $|G(\Gamma_{1})|=2^{|Q_{1}|}$ and 
$|G(\Gamma_{2})|=2^{|Q_{2}|}$.

On the other hand, by the uniqueness of the decomposition of a 
permutation into a product of disjoint cycles, a permutation
belongs to $G(\Gamma_{1})\cap G(\Gamma_{2})$ if and only if it is a 
product of transpositions belonging to both $G(\Gamma_{1})$ and 
$G(\Gamma_{2})$. Therefore, 
$$
G(\Gamma_{1})\cap G(\Gamma_{2})=\langle T(\Gamma_{1})\cap 
T(\Gamma_{2})\rangle=
\langle (i,j)\mid i\!\cdot \! j\in Q_{1}\cap Q_{2}\rangle,
$$
and then, arguing as in the previous paragraph, we see that 
$|G(\Gamma_{1})\cap G(\Gamma_{2})|=2^{|Q_{1}\cap Q_{2}|}$.

Now, it is well known that 
$$
|G(\Gamma_{1})\cdot G(\Gamma_{2})|=\frac{|G(\Gamma_{1})|\cdot 
|G(\Gamma_{2})|}{|G(\Gamma_{1})\cap G(\Gamma_{2})|},
$$
and hence
$$
d_{sgr}(\Gamma_{1},\Gamma_{2})= 
\ln \left(\frac{|G(\Gamma_{1})|\cdot 
|G(\Gamma_{2})|}{|G(\Gamma_{1})\cap G(\Gamma_{2})|^2}\right)=
\ln 2^{|Q_{1}|+|Q_{2}|-2|Q_{1}\cap Q_{2}|}=\ln 2^{|Q_{1}\Delta 
Q_{2}|},
$$
as we claimed.
\end{proof}

Notice in particular that, should Reidys and Stadler had defined 
their subgroup metric as 
$\log_{2}(|G(\Gamma_{1})\cdot G(\Gamma_{2})|/|G(\Gamma_{1})\cap 
G(\Gamma_{2})|)$, it would coincide with $|Q_{1}\Delta Q_{2}|$.

The third metric on $\SE$ proposed by Reidys and Stadler is actually a 
general way of defining metrics, rather than a single one, and it uses 
Magarshak and coworkers' algebraic representation of RNA secondary structures 
\cite{KMM93,Mag93,MB92}, recently extended in \cite{CMR} to cope with 
contacts other than Watson-Crick complementary base pairs.  These 
authors represent an RNA secondary structure $\Gamma=([n],Q)$ as an 
$n\times n$ complex symmetric matrix 
$S_{\Gamma}=(s_{i,j})_{i,j=1,\ldots,n}$
where
$$
s_{i,j}=\left\{
\begin{array}{rl}
-1 & \mbox{ if $i\neq j$ and $i\!\cdot\! j\in Q$}\\
1 & \mbox{ if $i=j$ and  $i\!\cdot\! l\notin Q$ for every $l$}\\
0 & \mbox{ otherwise}
\end{array}
\right.
$$
Since $S_{\Gamma}^{-1}=S_{\Gamma}$ for every $\Gamma\in \SE$, one can 
define for any  $\Gamma_{1},\Gamma_{2}\in \SE$ the \emph{transfer 
matrix} $T_{\Gamma_{1},\Gamma_{2}}=S_{\Gamma_{2}}\circ 
S_{\Gamma_{1}}$.  Then, Reidys and Stadler propose to measure the 
difference between two RNA secondary structures by defining a metric 
through
$$
(\Gamma_{1},\Gamma_{2})\mapsto \|T_{\Gamma_{1},\Gamma_{2}}\|,
$$
where $\|\cdot \|$ stands for some \emph{length function} on the group 
$GL(n,\CC)$ of $n\times n$ invertible complex matrices \cite[Def.\ 
9, Lem.\ 6]{RS96} (actually, Reidys and Stadler propose to use a matrix norm 
$\|\cdot \|$, but it is probably a misprint, as it would not yield a 
metric).  A simple and well-known length function on $GL(n,\CC)$ is 
$$
\|A\|=\mathrm{rank}(A-\mathrm{Id}),
$$
which allows to define a metric on $\SE$
$$
d_{mag}(\Gamma_{1},\Gamma_{2})=\mathrm{rank}(T_{\Gamma_{1},\Gamma_{2}}-\mathrm{Id}).
$$
This metric turns out to be equal to the involution metric $d_{inv}$ defined 
above.

\begin{proposition}
For every $\Gamma_{1},\Gamma_{2}\in \SE$,
$d_{mag}(\Gamma_{1},\Gamma_{2})=d_{inv}(\Gamma_{1},\Gamma_{2})$.
\end{proposition}

The proof of this proposition is similar to (and simpler than) the 
proof of \cite[Thm.\ 17]{CMR}, which establishes essentially this 
equality for the generalized algebraic representation of RNA secondary 
structures in the sense of Magarshak introduced in that paper, and 
therefore we omit it.

\end{document}